\documentclass[twoside, 11pt]{article}
\usepackage{latexsym}
\usepackage{amsmath}
\usepackage{amsfonts}
\usepackage{amssymb}
\renewcommand{\a }{\alpha }
 
\renewcommand{\d}{\delta }
\newcommand{\D }{\Delta }

\newcommand{\e }{\varepsilon }
\newcommand{\g }{\gamma}

\renewcommand{\l }{\lambda }

 \newcommand{\norm}[1]{\Vert {#1} \Vert} % norma
 \newcommand{\parent}[3]{\left #1 {#3} \right #2} % racchiude il testo tra il
 \newcommand{\barre}[1]{\parent \Vert \Vert {#1}} % racchiude il testo tra barre adeguate
\newcommand{\bara}[1]{\parent | | {#1}}
 \newcommand{\tonde}[1]{\parent (){#1}} % racchiude il testo tra tonde adeguate
 \newcommand{\quadre}[1]{\parent []{#1}} % racchiude il testo tra quadre adeguate
 \newcommand{\graffe}[1]{\parent \{ \}{#1}} % racchiude il testo tra graffe adeg.

\newcommand{\salt}{\noalign{\vskip .2truecm}}
\newcommand{\var }{\varphi }

\newcommand{\s }{\sigma }
\newcommand{\Sig }{\Sigma}

\renewcommand{\o }{\omega }
\renewcommand{\O }{\Omega }

\newcommand{\ov}{\overline}

\newcommand{\be}{\begin{equation}}
\newcommand{\ee}{\end{equation}}
\newenvironment{pf}{\noindent{\sc Proof}.\enspace}{\rule{2mm}{2mm}\medskip}

\newcommand{\R}{\mathbb{R}}

\newcommand{\Z}{\mathbb{Z}}

\newcommand{\N}{\mathbb{N}}
\newtheorem{lem}{Lemma}[section]
\newtheorem{pro}[lem]{Proposition}
\newtheorem{thm}[lem]{Theorem}
\newtheorem{rem}[lem]{Remark}

\begin{document}
\baselineskip3.3ex
%\author{{\sc Fran\c{c}ois ebobisse and Mohameden Ould Ahmedou}
%\thanks{SISSA, via Beirut
%2-4, Trieste 34014 Italy. Email: ebobisse@sissa.it\qquad ahmedou@sissa.it. The two authors are supported by  postdoctoral fellowships from SISSA.}}

\title{\bf On a nonlinear fourth order elliptic equation\\ involving
 the critical Sobolev exponent}

\author{
{\sc Fran\c{c}ois Ebobisse} \\
SISSA, Via Beirut 2-4\\ 
 34014 Trieste, Italy \\
{\it E-mail: ebobisse@sissa.it} \\ \\
{\large and}\\ \\
{\sc Mohameden Ould Ahmedou}\\ 
Rheinische Friedrich-Wilhelms-Universit\"{a}t Bonn \\
Mathematisches Institut \\
Beringstrasse 4, D-53115 Bonn\\
{\it E-mail: ahmedou@math.uni-bonn.de}}

\date{}
\maketitle

%{\sc Fran\c{c}ois ebobisse \qquad\qquad Mohameden Ould Ahmedou } \and { }}
%\begin{document}

%\maketitle

{\footnotesize
\begin{abstract}

\noindent 
In this paper we use an algebraic topological argument due to Bahri and Coron  to show how the topology of the domain influences the 
existence of positive solutions of some fourth order elliptic equation involving the critical Sobolev exponent.

\bigskip\bigskip


\noindent  2000 (AMS) Mathematics Subject Classification: 35J60, 35J65, 58E05  	
\end{abstract}
}
%\footnotetext{
%Address: SISSA, Via Beirut 2-4, 34014 Trieste, Italy, 
%E-mail: ebobisse@sissa.it \qquad\qquad ahmedou@sissa.it}
\section{Introduction}
We  are looking for the solution of the following problem
$$\left\{
\begin{array}{cc}
%\bega\label{pro1}
 \Delta^2u=u^{\frac{n+4}{n-4}}&\mbox{ in }\quad\O\\
\salt
 u>0& \mbox{ in }\quad\O\\
\salt
 u=\Delta u=0&\,\,\,\,\mbox{ on }\quad\,\partial\O
\end{array}
\right.\leqno (P)
$$
where $\O $ is a bounded domain of $\R ^n$ ($n\geq 5$) with a smooth boundary $\partial\O$.\\
The interest in this equation grew up from its ressemblance to some geometric equations involving Paneitz operator and widely studied in these last years (see for instance \cite{cy}, \cite{bcy}, \cite{csya}, \cite{cgy}, \cite{dhl}, \cite{nsyj}).
In contrast with the subcritical case, the variational problem corresponding to $(P)$ presents a lack of compactness which can be easily seen in the fact that the Euler functional associated to $(P)$ does not satisfy Palais-Smale condition (P.S.). This means that there exist non compact sequences along which the functional is bounded and its gradient goes to zero. This is due to the non compactness of the embedding $H^2(\O)\cap H^1_0(\O)\hookrightarrow L^{\frac{2n}{n-4}}(\O)$.\\
A first result concerning $(P)$ was obtained by Van der Vorst \cite{vander2}, who proved that $(P)$ has no solution when $\O$ is starshaped. Such non existence results were proved by Pohozaev \cite{poho}, for Yamabe type problems on a domain of $\R^n$. Pushing further the ressemblance of the two problems, the Yamabe type  from one part and the Paneitz type from the other part, we are led to investigate the influence of the topology of the domain $\O$ on the existence of the solution. We obtain a similar result as the one obtained by Bahri-Coron \cite{baco} for Yamabe type problems. Namely, if $k$ is a positive integer and we denote by $H_k(\O, {\Z}_2)$ the homology of dimension $k$ of $\O $ with ${\Z}_2$-coefficients,  we have the following result
\begin{thm}\label{main}
If $H_k(\O, {\Z}_2)\neq 0$ for some positive integer $k$, then the problem ${\rm (}P{\rm )}$ has a solution.
\end{thm}
The solution obtained in Theorem \ref{main} can be characterized to be a solution of higher energy and higher index. Precisely, we have the following characterization by Chen \cite{yc} of the solution obtained by the Bahri-Coron topological scheme.
\begin{thm}\label{charact}
The solution $u$ obtained in Theorem \ref{main} satisfies, for some positive integer $p_0$:
\begin{itemize}
\item[(i)] $p_0^{\frac{4}{n-4}}S\,\leq\,J(u)\,\leq\,(p_0+1)^{\frac{4}{n-4}}S;$
\item[(ii)] $ind\,(J,u)\,\leq\,(p_0+1)k\,+\,p_0,\quad ind\,(J,u)\,+\,dim\,ker\,\partial ^2J(u)\,\geq\,p_0k\,+\,p_0;$
\item[(iii)] $u$ induces some difference of topology at the level $p_0^{\frac{4}{n-4}}S$ 
where $S^{-1}$ is the best constant of the Sobolev embedding $H^2(\R^n)\hookrightarrow L^{\frac{2n}{n-4}}(\R ^n)$.
\end{itemize}
\end{thm} 
The proof of Theorem \ref{main} goes along the topological arguments introduced by Bahri-Coron \cite{baco}. That is, arguing by contradiction, we assume throughout the paper that the problem ($P$) has no solution. Then we introduce a family of continuous maps of pairs which induces  nontri\-vial homomorphisms in the relative homology under the assumption that ($P$) has no solution. Expanding carefully the Euler functional associated to ($P$), near its potential critical points at infinity, we prove through a study of the interaction between the solutions of the problem ``at infinity'' that, when some parameters of the above family become very large, the   homomorphisms induced become trivial, leading then to a contradiction. \\
It is worthwile to point out that our solution is obtained exploiting the topology of the domain, however the expansion of the functional involves the Green's functions for the bilaplacian operator under Navier conditions which suggests that the geometry of the domain would play some role in a necessary and sufficient condition to obtain a solution of ($P$).  Such feature is shared by the corresponding Yamabe type problem see \cite{blr}.\\
%We plan to prove this theorem by contradiction, using a topological argument introduced by Bahri-Coron \cite{baco}. That is, we will assume througout the paper that the problem ($P$) has no solution. Then we introduce a functional whose positive critical points in the unit sphere of some suitable Banach space are related to the solutions of the problem ($P$). Putting toghether the theory developped by Bahri \cite{ba} to find critical points at infinity, with some topological arguments in \cite{baco}, we prove that the assumption that ($P$) has no solution may lead to a topological contradiction. This paper is organized as follows........
The plan of the paper is the following. In section 2 we present some  technical tools, while in section 3 we study the expansion of the Euler functional associated to ($P$), near its potential critical points at infinity. Section 4 is devoted to the proof of Theorem \ref{main}. Finally in Section 5 we collect some useful lemmas.
%\centerline{\bf Acknowledgements:}
%\vskip .1truecm
%\noindent The authors would like to thank S.I.S.S.A for its support through po%stdoctoral fellowships.
%The research of the two authors is supported by S.I.S.S.A postdoctoral fellows%hips. 
 \section{Notation and Preliminaries}\label{Noprel}
Let $\O $ be a bounded domain of $\R ^n$ ($n\geq 5$) with a smooth boundary $\partial\O$. Let $E:=H^1_0(\O)\cap H^2(\O)$ equipped with the norm 
$$\norm{u}_{E}:=\norm{\Delta u}_{L^2(\O )}.$$
Let $$\Sigma :=\{u\in E:\, \norm{u}_{E}=1\}\quad\mbox{and}\quad \Sigma _+:=\{u\in \Sigma:\, u>0\}.$$
 Instead of working with the functional
$$J_0(u):={1\over 2}\displaystyle\int_\O (\Delta u)^2dx-{n-4\over 2n}\int_\O (u^+)^{2n\over n-4}dx$$
whose critical points are solutions of ($P$), it's more convenient here to work with the functional
\be\label{funct}
J(u):=\displaystyle{\tonde{\displaystyle\int_\O (\Delta u)^2dx}^{n\over n-4}\over 
\displaystyle\int_\O u^{2n\over n-4}dx}.
\end{equation}
One can easily verify that if $u$ is a critical point of $J$ on $\Sigma _+$, then $J(u)^{n-4\over 8}u$ is a solution of ($P$). Let us introduce the gradient flow of the functional $J$ which will be used in the classical deformation argument.
\begin{eqnarray}\label{gradf}
\begin{cases}
\displaystyle{\partial\eta(t,u)\over\partial t}=-J'(\eta(t,u)) &\quad (t,u)\in [0,+\infty[\times E\\
\salt
 \eta(0,u)=u & \quad\mbox{ }
\end{cases}
\end{eqnarray}
By the maximum principle, the above flow preserves the positive cone, namely we have
\begin{lem}\label{invflow}
$\Sigma _+$ is invariant with respect to $\eta(t,\cdot)$, i.e., if $u\in\Sigma _+$ then $\eta(t,u)\in\Sigma _+$.
\end{lem}
The functional $J$ is known to  not satisfy the Palais-Smale condition (PS for short) on $\Sigma _+$. However, as in Struwe \cite{struwe0} (see also \cite{struwe}) sequences of Palais-Smale failing the PS condition can be described. To this end, we introduce the solutions of the problem
$$
\begin{cases}
 \Delta^2v=v^{{n+4}\over {n-4}}&\mbox{ in }\quad \R^n \\
\salt
 v>0& \mbox{ in }\quad \R^n \\
\end{cases}
$$
which, according to Chang-Shou Lin \cite{CSL} are radially symmetric about some point $a\in\R^n$ and are given by the following family of $(n + 1)$ parameters
$$\delta_{a,\lambda}(x):=c_n\tonde{\lambda\over 1+\lambda ^2|x-a|^2}^{n-4\over 2}$$
where $c_n$ is a normalization factor depending only on $n$ and is such that $\int_{\R^n}\bigl(\Delta\delta_{a,\lambda}\bigr)^2dx=1$.\\
For any $a\in\O$, $\l >0$, Let $P\delta _{a,\l }$ be the unique solution of 
$$
\begin{cases}
 \Delta^2P\delta _{a,\l }=\delta _{a,\l }^{{n+4}\over {n-4}}&\mbox{ in }\quad \O \\
\salt
 P\delta _{a,\l }=\Delta P\delta _{a,\l }=0 & \mbox{ on }\quad \partial\O 
\end{cases}
$$
and set $\var _{a,\l}:= \delta _{a,\l }- P\delta _{a,\l }$. \\
Given a positive integer $p$  and $\e >0$, we denote by
$$
V(p,\e) :  =    \left\{ \begin{array}{cc}
   u \in \Sigma_+   \mbox{ such that } \exists \,  (a_1,\cdots,a_p) \in \O^p \mbox{ and  } \exists \, (\l_1,\cdots,\l_p) \in (\R^*_+)^p \mbox{ such that }   & \hskip -.4truecm\\
\salt
\barre{u - \frac{\sum_{i = 1}^p P\delta _{a_i,\l_i}}{ ||\sum_{i = 1}^pP\delta _{a_i,\l_i}  ||}} < \e , \, \mbox{ with } \, \l_i {\rm dist}\,(a_i,\partial\O) \geq \frac{1}{\e} \, \mbox{ and } \e_{ij} < \e \, &  \hskip -.4truecm
\end{array}
\right\}
$$
where $$ \e_{ij} \, = \, \frac{1}{\frac{\l_i}{\l_j} + \frac{\l_j}{\l_i} + \l_i \l_j|a_i - a_j|^2} $$
$dist$ being the euclidean distance in $\R^n$.
If $u$ is a function in $V(p,\e) $, one can find an optimal representation, arguing as in Proposition 7  of \cite{baco} , namely  we have:
\begin{lem}\label{l:rep}
For every $p \geq 1$, there exists $ \e > 0$ such that for any $ u \in V(p,\e)$ the minimization problem
$$
\inf_{\a_i,b_i,\mu_i} \, \barre{u - \sum_{i = 1}^p \a_i P\delta_{b_i,\mu_i}}
$$
has a unique solution, up to a permutation on the set of indices $ \{1,\cdots,p\}$, where $\a_i>0$, $b_i\in \overline\O$, $\mu_i>0$ for any $i=1,\cdots,p$.
\end{lem}
Setting $ b_p: = (p)^{\frac{4}{n - 4}}S,$ we have the following characterization of Palais-Smale sequences.
%S:=\tonde{\int_{\R^n}\delta _{a,\l }^{{2n}\over {n-4}}dx}^{4\over {n-4}}\quad\mbox{ and }\quad 
 %b_p: = (p)^{\frac{4}{n - 4}}S,$
%one can easily verify that $$S^{n-4\over 4} = c_n^{{2n}\over {n-4}}\int_{\R^n }\, \frac{1}{(1 + |x|^2)^n}dx.$$ 
\begin{pro}\label{p:ps}
Under the assumption that (P) has no solution, 
let $ (u_k) \,  \subset \Sig_+$ be a sequence satisfying $ J(u_k) \to c $, a positive number and $ \partial J (u_k) \to 0$. There exist an integer $ p \geq 1$ and a sequence $ (\e_k)_k$ such that $ u_k \in V(p,\e)$.
Conversely, let $ p \in \N^+$, let $ (\e_k)$ be a positive sequence with $ \lim_{k \to + \infty} \e_k = 0$ and let $ (u_k)\,\subset V(p,\e)$ then $\partial J(u_k) \to 0$ and $ J(u_k) \to b_p$.
\end{pro}
\begin{rem}\label{r:ps}
{\rm As observed by Bahri-Coron in \cite{baco}, this kind of conclusion were already obtained in other situations. For instance Sacks-Uhlenbeck \cite{SACUH}, Meeks-Yau \cite{MEYA} for Harmonic maps, Wente \cite{WENT}, Brezis-Coron \cite{breco} for $H$-systems.}
\end{rem}
Setting
$$
 W_p: = \{ u \in \Sig _+ , \mbox{ such that } J(u) < b_{p + 1} \}
$$
we state the following deformation lemma, whose proof is similar to Lemma 17 in Bahri-Brezis \cite{bb}.
\begin{lem}\label{defor}
Given $p\in\N^*$ and  $ \e > 0$. 
Under the asumption that (P) has no solution, the pair $ (W_p,W_{p - 1})$ retracts by deformation onto the pair $ (W_{p - 1} \cup A_p, W_{p - 1}) $ where $ A_p \subset V(p,\e) $.
\end{lem}

\section{Expansion of the functional near its potential critical point at infinity}\label{expan}

We fix $K$ a compact subset in $\O$. Let us recall that $O(x)$ will denote
functions such that $O(x)\leq C|x|$ for some postive constant $C$. 
%\begin{pro}\label{exp:prop} There exists an
%integer $ p_0$ and a positive real number $ \l_0 > 0$ such that  for any $
%(\a_1, \cdots, \a_{p_0})$ satisfying $ \a_i \geq 0, \, \sum_{i = 1}^{p_0} \,
%\a_i = 1$, for any $ (a_1,\cdots,a_{p_0}) \in K^p_0$, for any $ \l \geq \l_0$,
%we have $$ J\Bigl(\sum_{i = 1}^{p_0} \, \a_i \,P\delta_{\a_i, \l}\Bigr) \,
%\leq \, (p_0S)^{\frac{4}{n - 4}}  $$ \end{pro}
%To prove this propositon we will need the following lemmas.
This section being devoted to an expansion of $J$ in $V(p,\e)$, a first estimate is given in the following lemma.
\begin{lem}\label{lemB1}
For every $p \in \N^*$ and every $ \e > 0$, there exists $ \l_p = \l(p, \e )$ such that for any $ (\a_1, \cdots,\a_p)$ satisfying  $ \a_i \geq 0, \sum_{i = 1}^p \, \a_i \, = 1$, for any $ (a_1, \cdots,a_p) \in \,K^p$, for any $ \l \geq \l_p$, we have:
$$
J\Bigl(\sum_{i = 1}^{p} \, \a_i\, P\delta_{\a_i, \l}\Bigr) \, \leq \, \bigl(p+\e \bigr)^{\frac{4}{n - 4}}S. 
$$

\end{lem}
\begin{pf}
By adapting the lemma B2 of \cite{baco} in our situation, we get
$$
J\Bigl(\sum_{i = 1}^p \,  \a_i\, P\delta_{\a_i, \l}\Bigr) \, \leq \,\graffe{\frac{\displaystyle\int_\O \, \Bigl(\sum_{i = 1}^p\, \a_i\,\d_{a_i,\l}\Bigr)^{\frac{2n}{n - 4}}dx}{\displaystyle\int_{\O}\Bigl(\sum_{i = 1}^p \,\a_i\,P\d_{a_i,\l}\Bigr)^{\frac{2n}{n - 4}}dx } }^{\frac{1}{2}} \,\left[\int_\O \, \Bigl(\sum_{i = 1}^p\, \beta _i\,\d_{a_i,\l}\Bigr)^{\frac{2n}{n - 4}}dx \right]^{\frac{4}{n-4}}
$$
where $\beta _i:=\frac{\a_i\d_{a_i,\l}}{\sum_{j = 1}^p\, \a_j\d_{a_j,\l}}$. Now, from the convexity of $x\to |x|^{\frac{2n}{n - 4}}$ we have that

$$ \Bigl(\sum_{i = 1}^p\, \beta _i\,\d_{a_i,\l}\Bigr)^{\frac{2n}{n - 4}}\leq \sum_{i = 1}^p\, \beta _i\,\d_{a_i,\l}^{\frac{2n}{n - 4}}.$$
So, by integrating and using the definition of $S$ we obtain that
\begin{eqnarray*}
\int_\O \, \Bigl(\sum_{i = 1}^p\, \beta _i\,\d_{a_i,\l}\Bigr)^{\frac{2n}{n - 4}}\,dx &\leq &\int_\O \, \sum_{i = 1}^p\, \beta _i\,\d_{a_i,\l}^{\frac{2n}{n - 4}}\,dx\\
\salt
&\leq & pS^{n-4\over 4}\,+\, \int_{\R^n\setminus\O} \, \sum_{i = 1}^p\, \beta _i\,\d_{a_i,\l}^{\frac{2n}{n - 4}}\,dx \\
\salt
&\leq & pS^{n-4\over 4}\,+\, O\Bigl(\frac{1}{\l^n}\Bigr).
\end{eqnarray*}
Now using Lemma \ref{a1}, the inequality in Lemma \ref{l:ineq} (i) and the
convexity of $x\to |x|^{\frac{n+4}{n - 4}}$ we also get $$
J\Bigl(\sum_{i = 1}^p \,  \a_i\, P\delta_{\a_i, \l}\Bigr)\,\leq\, \graffe{\displaystyle\int_\O \, \Bigl(\sum_{i = 1}^p\, \a_i\,\d_{a_i,\l}\Bigr)^{\frac{2n}{n - 4}}dx\over \displaystyle\int_\O \, \Bigl(\sum_{i = 1}^p\, \a_i\,\d_{a_i,\l}\Bigr)^{\frac{2n}{n - 4}}dx}^{\frac{1}{2}}{p^{4\over n-4}S\Bigl(1\,+\, O\bigl(\frac{1}{\l ^{n-2}}\bigr)\Bigr)\over \Bigl(1\,+\, O\bigl(\frac{1}{\l ^{n-4}}\bigr)\Bigr)}
.$$
Then, for $\l\geq\l_p$, it follows that 
 $$J\Bigl(\sum_{i = 1}^p \,  \a_i\, P\delta_{\a_i, \l}\Bigr)\,\leq\,p^{4\over n-4}S\Bigl(1\,+\, O\bigl(\frac{1}{\l ^{n-2}}\bigr)\Bigr)\,\leq\,(p+\e)^{4\over n-4}S.$$
\end{pf}

In the sequel, we will prove that $\e$ can be taken equal to zero in the above lemma provided that $p$ is large.
\begin{lem}\label{b2}
For any integer $p\in [2,+\infty[$, there exist $\e_1>0$, $\l_2\in (0,\infty)$ such that for any $(a_1,\cdots,a_p)\in K^p$, $\l\in (\l_2,\infty)^p$, $(\a_1,\cdots,\a_p)\in (0,\infty)^p$, with $\sum_{i=1}^p\a_i=1$, if there exists $i_0\in\{1,\cdots,p\}$ such that $\a_{i_0}\leq\e_1$,
then $J\Bigl(\sum_{i = 1}^p \,  \a_i\, P\delta_{a_i, \l}\Bigr)\,\leq\,p^{4\over n-4}S$.
\end{lem}
\begin{pf}
Assume that $\a_p\neq 1$ and set $\tilde\a_i:=\frac{\a_i}{\sum_{i=1}^{p-1}\a_i}$, then by easy computations, one can see that
$$
N^{n-4\over n}:=\int_\O \, \Bigl(\sum_{i = 1}^p\, \a_i\,\Delta\d_{a_i,\l}\Bigr)^2dx = (1-\a_p)^2\int_\O \, \Bigl(\sum_{i = 2}^{p-1}\,\tilde \a_i\,\Delta\d_{a_i,\l}\Bigr)^2dx + O\Bigl(\a_p\Bigr)
$$
and
$$
D:=\int_\O \, \Bigl(\sum_{i = 1}^p\, \a_i\,\d_{a_i,\l}\Bigr)^{\frac{2n}{n - 4}}dx=(1-\a_p)^{2n\over n-4}\int_\O \, \Bigl(\sum_{i = 1}^p\, \tilde\a_i\,\d_{a_i,\l}\Bigr)^{\frac{2n}{n - 4}}dx +O\Bigl(\a_p^{\frac{2n}{n - 4}}\Bigr).
$$
Therefore
$$
J\Bigl(\sum_{i = 1}^p \,  \a_i\, P\delta_{a_i, \l}\Bigr)\,\leq\,J\Bigl(\sum_{i = 1}^{p-1} \,  \tilde\a_i\, P\delta_{a_i, \l}\Bigr)\Bigl(1+O\Bigl(\a_p^{\frac{n}{n - 4}}\Bigr)\Bigr).$$
For any $\e>0$, the exists $\e_1>0$ such that if $\a_p\,\leq\,\e_1$, then
$$
J\Bigl(\sum_{i = 1}^p \,  \a_i\, P\delta_{a_i, \l}\Bigr)\,\leq\,J\Bigl(\sum_{i = 1}^{p-1} \,  \tilde\a_i\, P\delta_{a_i, \l}\Bigr)(1+\e).$$
Now from Lemma \ref{lemB1} we have that
$$
J\Bigl(\sum_{i = 1}^p \,  \a_i\, P\delta_{a_i, \l}\Bigr)\,\leq\,\bigl(p-1 +\e\bigr)^{4\over n-4}S +\e\,\leq\,p^{4\over n-4}S.$$
\end{pf}

Let us recall that in the last lemma, we used the assumption that at least one $\a_i$ is small enough.  
Following \cite{baco}, in order to get an estimate of $J\Bigl(\sum_{i = 1}^p \,  \a_i\, P\delta_{a_i, \l}\Bigr)$ similar to the one in Lemma \ref{b2}, when all the $\a_i$ are bounded from below, we will need an expansion of $J\Bigl(\sum_{i = 1}^p \,  \a_i\, P\delta_{a_i, \l}\Bigr)$ involving the Green's function and its regular part for the bilaplacian operator under Navier boundary conditions (see Appendix).\\
% when $\l\min_{i\neq j}|a_i-a_j|$ is large enough.\\
We have the following proposition.
\begin{pro}\label{expagreen}
There exists a positive constant $C(p)$ such that
%\be\label{expagreen1}
$$
\bara{J\Bigl(\sum_{i = 1}^p \,  \a_i\, P\delta_{a_i, \l}\Bigr)\,-\,\Psi(\a,a,\l)}\,\leq\,\frac{C(p)}{\l^{n-2}d^{n-2}},
$$
%\ee
where $d=d(a):=\min_{i\neq j}|a_i-a_j|$,
\begin{eqnarray}\label{Psi}
%\label{expagreen2}
& { }&\hskip -.1truecm \Psi(\a,a,\l):\,= \,\frac{S|\a|^{2n\over n-4}}{||\a||^{2n\over n-4}}\left\{1\,-\,\frac{c_1}{\l^{n-4}}\left [\sum_{i=1}^pH(a_i,a_i)\left(\frac{\a_i^2}{|\a|^2}\,
-\,\frac{2\a_i^{2n\over n-4}}{||\a||^{2n\over n-4}}\right)\right.\right.\\
\salt
\nonumber &{ } &\hskip 6.8truecm  \left.\left. +\,\sum_{i\neq j}\tonde{\frac{2\a_i^{n+4\over n-4}\a_j}
{||\a||^{2n\over n-4}}\,-\,\frac{\a_i\a_j}{|\a|^2}}G(a_i,a_j)\right ]\right\},
\end{eqnarray}
with $G$ and $H$ the Green function and its regular part for  the bilaplacian operator under Navier boundary conditions (see Appendix), 
$$|\a|:=\Bigl(\sum_{i=1}^p\a_i^2\Bigr)^{1\over 2},\qquad 
||\a||:=\Bigl(\sum_{i=1}^p\a_i^{2n\over n-4}\Bigr)^{n-4\over 2n},\qquad
c_1:={nc_2c_n^{n+4\over n-4}\over  (n-4)S^{n-4\over 4}}$$
and $c_2$ is given by
\begin{equation}\label{c_2}
c_2:=\int_{\R^n}\frac{dy}{\bigl(1+|y|^2\bigr)^{n+4\over 2}}.
\end{equation}
\end{pro}
\begin{pf}
We set $u:=\displaystyle\sum_{i = 1}^p \,  \a_i\, P\delta_{a_i, \l}$ and we recall that $
J(u)=\displaystyle{\tonde{\displaystyle\int_\O (\Delta u)^2dx}^{n\over n-4}\over 
\displaystyle\int_\O u^{2n\over n-4}dx}.
$\\
Let us first begin with the expansion of the numerator. So,
\begin{eqnarray}\label{exp:N}
N^{n-4\over n}&=&\int_\O (\Delta u)^2\,dx\, =\, \int_\O\Delta^2 u\,u\,dx\\
\salt
\nonumber  & = &\sum_{i=1}^p\a_i^2\int_\O\delta_{a_i,\l}^{n+4\over
n-4}\,P\delta_{a_i,\l}\,dx\,+\, \sum_{i\neq
j}\a_i\a_j\int_\O\delta_{a_i,\l}^{n+4\over n-4}\,P\delta_{a_j,\l}\,dx.
\end{eqnarray}
Using the estimate in Lemma \ref{usefull} (ii), we get
\begin{equation}\label{exp:1}
\int_\O\delta_{a_i,\l}^{n+4\over
n-4}\,P\delta_{a_i,\l}\,dx\,=\,S^{n-4\over 4}\,-\,\frac{c_2c_n^{n+4\over
n-4}}{\l^{n-4}}H(a_i,a_i)\,+\,O\Bigl({1\over\l^{n-2}}\Bigr).
\end{equation}
On the other hand, given $i\neq j$,
\begin{eqnarray}\label{exp:2}
\int_\O\delta_{a_i,\l}^{n+4\over n-4}\,P\delta_{a_j,\l}\,dx & =
&\int_{\R^n}\delta_{a_i,\l}^{n+4\over n-4}\,P\delta_{a_j,\l}\,dx\,-\,
\int_{\R^n\setminus\O}\delta_{a_i,\l}^{n+4\over n-4}\,P\delta_{a_j,\l}\,dx\\
\salt
\nonumber & = &\int_{\R^n}\delta_{a_i,\l}^{n+4\over
n-4}\,P\delta_{a_j,\l}\,dx\,+\,O\bigl({1\over\l^n}\bigr)\\
\salt
\nonumber & = & \int_{\R^n}\delta_{a_i,\l}^{n+4\over
n-4}\,\delta_{a_j,\l}\,dx\,-\,   \int_{\R^n}\delta_{a_i,\l}^{n+4\over
n-4}\,\var_{a_j,\l}\,dx\,+\,  O\bigl({1\over\l^n}\bigr). 
\end{eqnarray}
We set $a_{ij}:=a_i\,-\,a_j$ and let us now estimate $I:=\int_{\R^n}\delta_{a_i,\l}^{n+4\over n-4}\,\delta_{a_j,\l}\,dx$. By
easy computations, we have 
\begin{eqnarray}\label{exp:I}
I&=&c_n^{2n\over
n-4}\int_{\R^n}\left(\frac{\l}{1\,+\,\l^2|x\,-\,a_i|^2}\right)^{n+4\over
2}\left(\frac{\l}{1\,+\,\l^2|x\,-\,a_j|^2}\right)^{n-4\over
2}\,dx\\ 
\salt
\nonumber &=&c_n^{2n\over
n-4}\int_{\R^n}\frac{1}{\bigl(1\,+\,|y|^2\bigr)^{n+4\over
2}\bigl(1\,+\,|y\,-\,\l a_{ij}|^2\bigr)^{n-4\over
2}}\,dy.
\end{eqnarray}
We have also
$$ 1\,+\,|y\,-\,\l a_{ij}|^2=\bigl(1\,+\,\l^2|
a_{ij}|^2\bigr)\left\{1+\frac{|y|^2\,-\,2\l y\cdot a_{ij}}{1\,+\,\l^2|
a_{ij}|^2}\right\}
$$
hence, for $\l d(a)$ large and $|y|\leq {1\over 4}\l |a_{ij}|$ we obtain
\begin{eqnarray}\label{exp:3}
&{ } &\hskip -2truecm \bigl(1\,+\,|y\,-\,\l
a_{ij}|^2\bigr)^{-{n-4\over 2}}=\bigl(1\,+\,\l^2|
a_{ij}|^2\bigr)^{-{n-4\over 2}}\left\{1+\frac{(n-4)\l y\cdot
a_{ij}}{1\,+\,\l^2| a_{ij}|^2}\right.\\ 
\salt
\nonumber &{ } &\hskip 5.7truecm \left. +\,\,O\left (\frac{|y|^2}{1\,+\,\l^2|
a_{ij}|^2}\right)\right\}.
\end{eqnarray}
Let
$$A(y):\,=\,\left(\frac{1}{1\,+\,|y|^2}\right)^{n+4\over
2}\left(\frac{1}{1\,+\,|y\,-\,\l a_{ij}|^2}\right)^{n-4\over
2}.$$
Then
\begin{eqnarray}\label{exp:4}
\hskip .4truecm \int_{|y|\leq\l|a_{ij}|/4}A(y)\,dy
&=&\frac{1}{\bigl(1\,+\,\l^2| a_{ij}|^2\bigr)^{n-4\over
2}}\,\left\{\int_{|y|\leq\l|a_{ij}|/4}\left(\frac{1}{1\,+\,|y|^2}\right)^{n+4\over 2}\,dy\right.\\ \salt
\nonumber & { } &\,\hskip 1.1truecm\left.+\,\frac{1}{\bigl(1\,+\,\l^2|
a_{ij}|^2\bigr)}\,O\left(\int_{|y|\leq\l|a_{ij}|/4}\frac{|y|^2}{\bigl(1\,+\,|y|^2\bigr)^{n+4\over 2}}\,dy\right)
\right\}
\end{eqnarray}
From the following identities
$$
\int_{|y|\leq\l |a_{ij}|/4}\frac{|y|^2}{\bigl(1\,+\,|y|^2\bigr)^{n+4\over
2}}\,dy\,=\,O\Bigl(\frac{1}{\l^2|a_{ij}|^2}\Bigr)$$
\vskip .2truecm
$$
\int_{|y|\leq\l|a_{ij}|/4}\frac{1}{\bigl(1\,+\,|y|^2\bigr)^{n+4\over
2}}\,dy\,=\,\int_{\R^n}\frac{1}{\bigl(1\,+\,|y|^2\bigr)^{n+4\over
2}}\,dy\,+\,O\Bigl(\frac{1}{\l^4|a_{ij}|^4}\Bigr)$$
\vskip .2truecm
$$
\frac{1}{\bigl(1\,+\,\l^2| a_{ij}|^2\bigr)^{n-4\over
2}}\,=\,\frac{1}{\l^{n-4}| a_{ij}|^{n-4}}\,+\,O\Bigl(\frac{1}{\l^{n-2}|
a_{ij}|^{n-2}}\Bigr)$$
we finally obtain
\begin{equation}\label{exp:5}
\int_{|y|\leq\l|a_{ij}|/4}A(y)\,dy\,=\,\frac{c_2}{\l^{n-4}| a_{ij}|^{n-4}}\,+\,O\Bigl(\frac{1}{\l^{n-2}|
a_{ij}|^{n-2}}\Bigr)
\end{equation}
with $c_2$ given in (\ref{c_2}). \\
Now, in order to estimate $\int_{|y|>\l|a_{ij}|/4}A(y)\,dy$ we introduce the sets
$$
B_1:=\left\{y\in\R^n\mbox{ : }|y\,-\,\l a_{ij}|\,\leq\,{\l |a_{ij}|\over 4}\right\}\quad\mbox{ and }\quad B_2:=\left\{y\in\R^n\mbox{ : }|y|\,\leq\,{\l |a_{ij}|\over 4}\right\}.$$
Then
\begin{equation}\label{exp:6}
\int_{\R^n\setminus (B_1\cup B_2)}A(y)\,dy \,\leq \,{C\over \l^{n-4}| a_{ij}|^{n-4}}\int_{\l |a_{ij}|}^\infty {r^{n-1}\over \bigl(1\,+\,r^2\bigr)^{n+4\over 2}}\,dr\,=\, O\Bigl(\frac{1}{\l^n|
a_{ij}|^n}\Bigr).
\end{equation}
On the other hand,
\begin{equation}\label{exp:7}
\int_{B_1}A(y)\,dy \,\leq \,{C\over \l^{n+4}| a_{ij}|^{n+4}}\int^{\l |a_{ij}|}_0 {r^{n-1}\over \bigl(1\,+\,r^2\bigr)^{n-4\over 2}}\,dr\,=\, O\Bigl(\frac{1}{\l^n|
a_{ij}|^n}\Bigr).
\end{equation}
Finally, from (\ref{exp:5}), (\ref{exp:6}) and (\ref{exp:7}), we get
\begin{eqnarray}\label{exp:8}
\int_{\R^n}A(y)\,dy &=& \int_{\R^n\setminus (B_1\cup B_2)}A(y)\,dy \,+\,\int_{B_1}A(y)\,dy \,+\,\int_{B_2}A(y)\,dy \\
\salt
\nonumber &=&\frac{c_2}{\l^{n-4}| a_{ij}|^{n-4}}\,+\,O\Bigl(\frac{1}{\l^{n-2}|
a_{ij}|^{n-2}}\Bigr).
\end{eqnarray}
From  (\ref{exp:I}) and (\ref{exp:8}), we obtain
\begin{equation}\label{exp:9}
\int_\O\delta_{a_i,\l}^{n+4\over n-4}\,\delta_{a_j,\l}\,dx \,=\,\frac{c_2}{\l^{n-4}| a_{ij}|^{n-4}}\,+\,O\Bigl(\frac{1}{\l^{n-2}|
a_{ij}|^{n-2}}\Bigr).
\end{equation}
So, from (\ref{exp:2}), (\ref{exp:9}) and Lemma \ref{usefull} (ii) we obtain that
\begin{eqnarray}\label{exp:10}
\int_\O\delta_{a_i,\l}^{n+4\over n-4}\,P\delta_{a_j,\l}\,dx & =
&\frac{c_2c_n^{n+4\over n-4}}{\l^{n-4}}\quadre{\frac{c_n}{|a_i-a_j|^{n-4}}\,-\,H(a_i,a_j)}\,+\,O\Bigl(\frac{1}{\l^{n-2}|
a_{ij}|^{n-2}}\Bigr)\\
\salt
\nonumber & = & \frac{c_2c_n^{n+4\over n-4}}{\l^{n-4}}G(a_i,a_j)\,+\,O\Bigl(\frac{1}{\l^{n-2}d^{n-2}}\Bigr).
\end{eqnarray}
Putting toghether (\ref{exp:N}), (\ref{exp:1}) and (\ref{exp:10}) we finally get that
\begin{eqnarray}
\nonumber\hskip -2.truecm N^{n-4\over n}&=&|\a|^2S^{n-4\over 4}\,-\,{c_2c_n^{n+4\over n-4}\over\l^{n-4}}\sum_{i=1}^p\a_i^2H(a_i,a_i)\\
\nonumber&{ }&\hskip 1.5truecm +\,{c_2c_n^{n+4\over n-4}\over\l^{n-4}}\sum_{i\neq j}\a_i\a_j\,G(a_i,a_j)\,+\,O\Bigl(\frac{1}{\l^{n-2}d^{n-2}}\Bigr)
\end{eqnarray}
from which it follows that
\begin{eqnarray}\label{exp:final:N}
N & = & |\a|^{2n\over n-4}S^{n\over 4}\left\{1\,-\,{c_2c_n^{n+4\over n-4}\over \l^{n-4}|\a|^2S^{n-4\over 4}}\left[\sum_{i=1}^p\a_i^2H(a_i,a_i)\, -\,\sum_{i\neq j}\a_i\a_j\,G(a_i,a_j)\right]\right.\\
\salt
\nonumber & { } & \left.\hskip 4.8truecm +\,\,\,O\Bigl(\frac{1}{\l^{n-2}d^{n-2}}\Bigr)\right\}^{n\over n-4}
\end{eqnarray}

Now, let us write the expansion of $D:=\displaystyle\int_\O u^{2n\over n-4}\,dx\,=\,\displaystyle\int_\O\left(\sum_{i=1}^p\a_i\,P\delta_{a_i,\l}\right)^{2n\over n-4}\,dx$ \\
for $\l\,d(a)\to\infty$. To this end, we set 
$$l:=\min_i{\rm dist}\,(a_i,\partial\O)\quad\mbox{ and }\quad B_i:\,=\,\left\{y\in\O\mbox{: }|y-a_i|\,<\,d'\right\}\quad\mbox{with}\quad d':=\min( d(a)/2,l).$$
We have
\begin{eqnarray*}
\int_{\O\setminus\bigl(\cup_{i=1}^pB_i\bigr)}\Bigl(\sum_{i=1}^p\a_i\,P\delta_{a_i,\l}\Bigr)^{2n\over n-4}dx&\leq &\int_{\O\setminus\cup_{i=1}^pB_i}\sum_{i=1}^p\a_i\,\bigl(P\delta_{a_i,\l}\bigr)^{2n\over n-4}\,dx\\
\salt
& = &\int_{\O\setminus\cup_{i=1}^pB_i}\sum_{i=1}^p\a_i\,\bigl(\delta_{a_i,\l}\,-\,\var_{a_i,\l}\bigr)^{2n\over n-4}\,dx.
\end{eqnarray*}
Recalling that $\norm{\var_{a,\l}}_{\infty,\O}\,\leq\,O\Bigl(\frac{1}{\l^{n-4\over 2}}\Bigr)$ (see Lemma \ref{a1}), we obtain
\begin{eqnarray}\label{exp:D1}
\int_{\O\setminus\bigl(\cup_{i=1}^pB_i\bigr)}\Bigl(\sum_{i=1}^p\a_i\,P\delta_{a_i,\l}\Bigr)^{2n\over n-4}dx&\leq &C(p)\sum_{i=1}^p\int_{\O\setminus\bigl(\cup_{i=1}^pB_i\bigr)}\delta_{a_i,\l}^{2n\over n-4}\,dx\,+\,O\Bigl(\frac{1}{\l^{n}}\Bigr)\\
\salt
\nonumber & \leq & C(p)\int_{r\geq\l\,d'}{r^{n-1}\over (1\,+\,r^2)^n}\,dr\,+\,O\Bigl(\frac{1}{\l^{n}}\Bigr)\\
\salt
\nonumber & = &\,O\Bigl(\frac{1}{\l^{n}{d'}^n}\Bigr).
\end{eqnarray}
In order to estimate $\int_{B_i} u^{2n\over n-4}\,dx$, we write
$$\sum_{j=1}^p\a_j\,P\delta_{a_j,\l}\,=\,\a_i\,\delta_{a_i,\l}\,+\,\sum_{j\neq i}\a_j\,P\delta_{a_j,\l}\,-\,\a_i\,\var_{a_i,\l}.$$
Using Lemma \ref{l:ineq} (i), we get for $q>2$
\begin{eqnarray*}\label{exp:D12}
\Bigl(\sum_{j=1}^p\a_j\,P\delta_{a_j,\l}\Bigr)^q &=&\bigl(\a_i\,\delta_{a_i,\l}\bigr)^q\,+\,q\bigl(\a_i\,\delta_{a_i,\l}\bigr)^{q-1}\Bigl(\sum_{j\neq i}\a_j\,P\delta_{a_j,\l}\,-\,\a_i\,\var_{a_i,\l}\Bigr)\,\,+\\
& { } & \hskip - .8truecm+\,\,O\left(\Bigl(\sum_{j\neq i}\a_j\,P\delta_{a_j,\l}\,-\,\a_i\,\var_{a_i,\l}\Bigr)^q \right.\\
& { } & \left.\hskip 1.2truecm +\,\bigl(\a_i\,\delta_{a_i,\l}\bigr)^{q-2}\inf\Bigl(\bigl(\a_i\,\delta_{a_i,\l}\bigr)^2\,,\,\Bigl(\sum_{j\neq i}\a_j\,P\delta_{a_j,\l}\,-\,\a_i\,\var_{a_i,\l}\Bigr)^2\Bigr)\right)
\end{eqnarray*}
Now, for $q:={2n\over n-4}$, we have on $ B_i$ that
$$\Bigl(\sum_{j\neq i}\a_j\,P\delta_{a_j,\l}\,-\,\a_i\,\var_{a_i,\l}\Bigr)^q\,=\,O\Bigl(\frac{1}{\l^n{d'}^{2n}}\Bigr)$$
$$\bigl(\a_i\,\delta_{a_i,\l}\bigr)^{q-2}\Bigl(\sum_{j\neq i}\a_j\,P\delta_{a_j,\l}\,-\,\a_i\,\var_{a_i,\l}\Bigr)^2\,=\,O\left(\frac{\delta_{a_i,\l}^{8\over n-4}}{\l^n{d'}^{2n}}\right).$$
Putting toghether these estimates we get
\begin{eqnarray*}\label{exp:D2}
\int_{B_i} u^{2n\over n-4}\,dx &=&\a_i^{2n\over n-4}\int_{B_i} \delta_{a_i,\l}^{2n\over n-4}\,dx\,+\,{2n\over n-4}\a_i^{n+4\over n-4}\int_{B_i} \delta_{a_i,\l}^{n+4\over n-4}\Bigl(\sum_{j\neq i}\a_j\,P\delta_{a_j,\l}\,-\,\a_i\,\var_{a_i,\l}\Bigr)\,dx\\
\salt
& { } &\hskip 3truecm +\,\,O\left(\frac{1}{\l^n{d'}^{2n}}\int_{B_i}\delta_{a_i,\l}^{8\over n-4}dx\right),
\end{eqnarray*}
Recall that
$$ 
\int_{B_i}\delta_{a_i,\l}^{8\over n-4}dx\,=\,O\Bigl(\frac{{d'}^{n-8}}{\l^4}\Bigr).$$
So, analogously to the computations done in (\ref{exp:1}), (\ref{exp:10}) and Lemma \ref{usefull} we obtain
\begin{eqnarray}\label{exp:D3}
\int_{B_i} u^{2n\over n-4}\,dx &=&\a_i^{2n\over n-4}\,S^{n-4\over 4}\,+\,{2n\over n-4}{c_2c_n^{n+4\over n-4}\a_i^{n+4\over n-4}\over\l^{n-4}}\sum_{j\neq i}\a_jG(a_i,a_j)\\
\salt
\nonumber & { } & -\,\,\,{2n\over n-4}{c_2c_n^{n+4\over n-4}\a_i^{2n\over n-4}\over\l^{n-4}}H(a_i,a_i)\,+\,O\Bigl({1\over\l^{n-2}{d'}^{n-2}}\Bigr).
\end{eqnarray}
Finally, from (\ref{exp:D1}) and (\ref{exp:D3}) we get
\begin{eqnarray}\label{exp:D4}
\hskip -1.5truecm D&=&\sum_{i=1}^p\left[\a_i^{2n\over n-4}\,S^{n-4\over 4}\,+\,{2n\over n-4}{c_2c_n^{n+4\over n-4}\a_i^{n+4\over n-4}\over\l^{n-4}}\sum_{j\neq i}\a_jG(a_i,a_j)\right.\\
\salt
\nonumber & { } & \hskip 1truecm \left.-\,\,\,{2n\over n-4}{c_2c_n^{n+4\over n-4}\a_i^{2n\over n-4}\over\l^{n-4}}H(a_i,a_i)\right]\,+\,O\Bigl({1\over\l^{n-2}{d'}^{n-2}}\Bigr)\\
\salt
\salt
\nonumber &=& ||\a||^{2n\over n-4}S^{n-4\over4}\left[1\,+\,{2n\over n-4}{c_2c_n^{n+4\over n-4}\over\l^{n-4}||\a||^{2n\over n-4}S^{n-4\over 4}}\sum_{j\neq i}\a_i^{n+4\over n-4}\a_jG(a_i,a_j)\right.\\
\salt
\nonumber & { } & \hskip 1truecm \left.-\,\,\,{2n\over n-4}{c_2c_n^{n+4\over n-4}\over\l^{n-4}||\a||^{2n\over n-4}S^{n-4\over 4}}\sum_{i=1}^p\a_i^{2n\over n-4}H(a_i,a_i)\right]\,+\,O\Bigl({1\over\l^{n-2}{d'}^{n-2}}\Bigr).
\end{eqnarray}
Then (\ref{exp:final:N}) and  (\ref{exp:D4}) give
\begin{eqnarray}\label{exp:J1}
\nonumber\qquad J(u) &=&\left({|\a|\over ||\a||}\right)^{2n\over n-4}S\left\{1\,-\,{n\over n-4}\,{c_2c_n^{n+4\over n-4}\over \l^{n-4}|\a|^2S^{n-4\over 4}}\left[\sum_{i=1}^p\a_i^2H(a_i,a_i)\, -\,\sum_{i\neq j}\a_i\a_j\,G(a_i,a_j)\right]\right.\\
\salt
\salt
\nonumber & { } & \left.\hskip 3.0truecm +\,\,\,O\left(\frac{1}{\l^{n-2}d^{n-2}}\right)\right\}\,\times\\
\salt
\nonumber & { } &\times\,\left\{\left[1\,-\,{2n\over n-4}{c_2c_n^{n+4\over n-4}\over\l^{n-4}||\a||^{2n\over n-4}S^{n-4\over 4}}\sum_{j\neq i}\a_i^{n+4\over n-4}\a_jG(a_i,a_j)\right.\right.\\
\salt
\nonumber & { } & \hskip 1truecm \left.\left.+\,\,{2n\over n-4}{c_2c_n^{n+4\over n-4}\over\l^{n-4}||\a||^{2n\over n-4}S^{n-4\over 4}}\sum_{i=1}^p\a_i^{2n\over n-4}H(a_i,a_i)\right]\,+\,O\Bigl({1\over\l^{n-2}{d'}^{n-2}}\Bigr)\right\}.
\end{eqnarray}
Then, by easy computation we obtain 
$$J(u)\,=\,\Psi(\a,a,\l)\,+\,O\left(\frac{1}{\l^{n-2}d^{n-2}}\right)$$
with $\Psi(\a,a,\l)$ given in (\ref{Psi}).
\end{pf}

\section{Proof of Theorem \ref{main}}

For the  proof of the Theorem \ref{main}, we introduce the following notations:\\
For any $ p \geq 1$ and $ \l > 0$, let 
$$
B_p = B_p(K) = \left\{ \sum_{i = 1}^p \a_i \d_{a_{i}} \, , \a_i \geq 0, \, \sum_{i = 1}^p \a_i = 1 , \, a_i \in K \, \right\}\quad\mbox{and}\quad  B_0 = B_0(K) = \emptyset 
$$
where $K$ is a compact subset in $\O$. 
Let also $ f_p(\l)$  denote the map from $ B_p(K)$ to $ \Sig_+$ defined by 
$$
f_p(\l) \Bigl(\sum_{i = 1}^p \a_i \, \d_{a_{i}}\Bigr) \, = \, \frac{\sum_{i = 1}^p \, P\d_{a_{i},\l_i}}{\barre{\sum_{i = 1}^p \, P\d_{a_{i},\l_i}}}
$$
Clearly we have $ B_{p - 1} \subset B_p$ and $ W_{p - 1} \subset W_p$. Moreover $ f_p(\l)$ enjoys the following properties:

\begin{pro}\label{p:tt}
The function $ f_p(\l)$ has the following properties:
\begin{description}

\item(i) For any integer $ p \geq 1$, there exists a real number $\l_p > 0 $ such that 
$$
f_p(\l) : \, B_p(K) \to W_p \, \mbox{ for any } \l \geq \l_p
$$

\item(ii) 
There exists an integer $ p_0 > 1$, such that for any integer $ p \geq p_0$ , and for any $ \l \geq \l_{p_{0}}$, the map of pairs $ f_p(\l) : \, (B_p,B_{p - 1}) \, \to (W_p, W_{p - 1})$  satisfies $ (f_p)_*(\l) \equiv 0$
where 
$$
 (f_p(\l))_* : \,H_* (B_p,B_{p - 1}) \, \to H_*(W_p, W_{p - 1})
$$
and $H_* $ is the $*$th homology group with $ \Z_{2} $ coefficients.
\end{description}
\end{pro}

\begin{pf}
(i) is a direct consequence of Lemma \ref{b2}.\\
(ii) is a consequence of Proposition \ref{expagreen}. Indeed, there exist $c_3$ and $c_4$ such that:
$$ H(y,y)\,\leq\,c_3\mbox{ for all } y\in K$$
$$G(y_1,y_2)\,\geq\,c_4\mbox{ for all }y_1,\,y_2\in K^2.$$
Hence from Proposition \ref{expagreen}, one easily deduces the following:\\
There exist two positive real numbers $\ov{\nu}$ and $\ov{\eta}$ such that, for any positive integer number $p$, there exists a constant $C(p)$ such that, for any $\l\in [1,\infty[$ and any $a\in K^p$ with $d(a)\neq 0$,
$$\max_{\D_{p-1}}\Psi(\a,a,\l)\,\leq\,p^{4\over n-4}\left[S\,+\, {2\over \l^{n-4}}(\ov{\nu}-p\ov{\eta})\right]\,+\, {C(p)\over (\l d(a))^{n-2}},$$
where $\D_{p-1}:=\{(\a_1,\cdots,\a_p)\mbox{, }\a_i\geq 0,\,\sum_{i=1}^p\a_i=1\}$. 
On the other hand, let us point out that there exist some $d_0$  positive and $\l_3\geq 1$ such that $\Psi (\a,a,\l)\,\leq\,p^{4\over n-4}S$, indeed
$$\lim_{|y_1-y_2|\to 0}G(y_1,y_2)=+\infty.$$
Suppose that $d(a)\,\geq\,d_0$. Now, to prove Proposition \ref{p:tt}, it is enough to choose $p_0$ such that 
$$\ov{\nu}-p_0\ov{\eta}\,<\,0.$$
\end{pf}

Let us point out that
$H_k(\O)\neq 0$ implies that there exists a $k$-dimensional compact connected ${\mathcal C}^\infty$ manifold $V$ without boundary and a continuous map $h:V\to\O$ such that if $w_k$ denotes the class of orientation of $V$, then $h_*(w_k)\neq 0$ and there exists a compact  ${\mathcal C}^\infty$ manifold $K$ such that $h(V)\,\subset\,K\,\subset\,\O$ (see Thom \cite{thom}).\\
Let 
$$ F_p = \{ (a_1,\cdots,a_p) \in V^p \quad \mbox{such that } \exists \, i \not= j \quad \mbox{with } a_i = a_j \}.$$
Let $ \s_p$ be the symmetric group of order p, which acts on $ F_p$, and let $T_p $ be a $ \s_p$- equivariant tubular neighborhood of $F_p$, in $V^p $ (The existence of a such neighborhood is derived in the book of G. Bredon \cite{br}, see also Appendix C in \cite{baco}).\\
From another part,  considering the topological pair $(B_p(V),B_{p - 1}(V)) $ we observe that $(B_p(V) \setminus B_{p - 1}(V)) $ can be described as $ (V_0^p)^* \times _{\s_p} (\D_p \setminus \partial \D_{p - 1}) $ where 
$$(V_0^p)^* = \{ (a_1,\cdots,a_p) \in (\partial V)^p \quad \mbox{such that } a_i \not= a_j , \forall i \not= j \}$$
We notice that $ (V_0^p)^*\times _{\s_p} (\D_p \setminus \partial \D_{p - 1})  $ is a noncompact manifold of dimension $ kp + p - 1$.\\
For $ 0 < \theta < 1$,  let $ \mathcal{M}_p: = \tilde{V_p} \times _{\s_p} \D_{p - 1}^{\theta}$, where $ \tilde{V_p}: = \ov{V^p \setminus T_p }$ and 
$$\D_{p - 1}^{\theta} := \left\{ (\a_1,\cdots,\a_p) \in \D_{p - 1} \quad \mbox{such that } \frac{\a_i}{\a_j} \in [1 - \theta, 1 + \theta ], \forall i,\, j \right\} .$$
$ \mathcal{M}_p$ is a manifold which can be seen as a subset of $ B_p(V)$, and the topological pair $(B_p(V),\mathcal{M}^c) $ retracts by deformation onto $ (B_p(V), B_{p - 1}(V))$, we thus have 
$$
H_*(B_p(V),B_{p - 1}(V)) = H_*(B_p(V), \mathcal{M}^c_p)
$$
Thus by excision we have
$$
H_*(B_p(V),B_{p - 1}(V)) = H_*(\mathcal{M},\partial  \mathcal{M}_p)
$$
Since any manifold is orientable modulo its boundary with $ \Z_2$ coefficients, we have a nonzero orientation class in $ H_{kp + p - 1}(B_p(V),B_{p - 1}(V))$ which we denote by $ \o_p$.\\
In contrast with Proposition \ref{p:tt}, we have the following Proposition:
\begin{pro}\label{p:ta}
Under the assumption that (P) has no solution, we have,
$$
 \mbox{for every } \, p \in \N^* \quad (f_p'(\l))_*(\o_p) \not\equiv 0.
$$
where $f_p'(\l):B_p(V)\,\stackrel{h}{\longrightarrow}\,B_p(h(V))\,\stackrel{f_p(\l)}{\longrightarrow}\,W_p$.
\end{pro}

\begin{pf} An abstract topological argument displayed in \cite{baco} , pp 260-265 , see also \cite{bb}, which extends virtually to our framework shows that:
$\mbox{If } \, (f'_1(\l)_* \not\equiv 0 \quad \mbox{then } (f'_p(\l))_* \not\equiv 0 \, \mbox{for every } p \geq 2.$ Since $ J_{S + \e}$ , for $ \e > 0$ small enough satisfies $ J_{S + \e} \subset V(1,\d)$ , where $ \d \to 0$ if $ \e \to 0$ , one can define using Lemma \ref{l:rep} a continuous map $ s : J_{S + \e} \to K$ which associates to $  u = \ov{\a} P\delta_{\ov{a},\ov{\l}} + v \in J_{S + \e}$ a point $ \ov{a} \in K$.
Here $ (\ov{\a},\ov{a},\ov{\l})$ is the unique solution (see Lemma \ref{l:rep}) of the minimization problem: 
$$ \min \{ \barre{u - \a P\delta_{a,\l}} , \, \a \geq 0, \l>0, a \in K \}.$$
Let $ r: W_1 \to J_{S + \e}$ denote the retraction by deformation of $ W_1$ onto $ J_{S + \e}$; the existence of a such retraction by deformation follows from the assumption that (P) has no solution from one part and from Lemma \ref{defor} from another part.
Let us observe that $ s \circ  r \circ f'_1(\l) = id_V$
hence $ (f'_1(\l))_*(\o_1) \not\equiv 0$ , where $ \o_1$ is the orientation class of $ V$.
Therefore the proof of Proposition \ref{p:ta} is reduced to the abstract topological argument of Bahri-Coron \cite{baco}.
\end{pf}

\vskip .2truecm
\noindent { \bf { Proof of Theorem \ref{main} completed:}}
Proposition \ref{p:ta} is in contradiction with Proposition \ref{p:tt}. Therefore (P) has a solution and Theorem \ref{main} is thereby established.
\section{Appendix}\label{Appen}
Let $G$ be the Green's function for the  bilaplacian, that is, given $x\in\O$,
$$
\begin{cases}
 \Delta^2G(x,\cdot)=c_n\delta_x &\mbox{ in }\quad \O \\ 
\salt
G(x,\cdot)=\Delta G(x,\cdot)=0& \mbox{ on }\quad \partial\O
\end{cases}
$$
where $c_n$ is the normalization constant in the definition of $\d_{a,\l}$. If we set $$H(x,y):=\frac{c_n}{|x-y|^{n-4}}-G(x,y)$$
then 
$$
\begin{cases}
 \Delta^2H(x,\cdot)=0 &\mbox{ in }\quad \O \\ 
\salt
H(x,\cdot)=\displaystyle\frac{c_n}{|x-y|^{n-4}}& \mbox{ on }\quad \partial\O\\
\salt
\Delta H(x,\cdot)=-2(n-4)\displaystyle\frac{c_n}{|x-y|^{n-4}}& \mbox{ on }\quad \partial\O
.\end{cases}
$$
The following lemma gives the estimate of $\var_{a,\l}$ in terms of the function $H(a,\cdot)$.
\begin{lem}\label{a1}
Let $\var_{a,\l}:=\d_{a,\l}-P\,\d_{a,\l}$, then 
$$\var_{a,\l}(y)=\frac{1}{\l^{n-4\over
2}}H(a,y)\,+\,O\Bigl(\frac{1}{\l^{n\over
2}}\Bigr)\quad\mbox{on}\quad\O .$$ 
\end{lem} 
\begin{pf}
From the definition of $P\,\d_{a,\l}$ we have
$$
\begin{cases}
 \Delta^2\var_{a,\l}=0&\mbox{ in }\quad \O \\
\salt
 \var_{a,\l}=\d_{a,\l} & \mbox{ on }\quad \partial\O\\ 
\salt
\Delta \var_{a,\l}=\Delta \d_{a,\l}& \mbox{ on }\quad \partial\O.
\end{cases}
$$
Now if we set $g_{a,\l}(y):=\displaystyle\frac{1}{\l^{n-4\over 2}}H(a,y),$
then by the maximum principle we obtain that 
$$\norm{\var_{a,\l}\,-\,g_{a,\l}}_{\infty,\O}\,\leq\,C\norm{\Delta\delta_{a,\l}\,-\,\Delta g_{a,\l}}_{\infty, \partial\O}=O\Bigl(\frac{1}{\l ^{n\over 2}}\Bigl)\quad\mbox{for some positive constant }C.$$ 
Moreover, we have the estimate $\norm{\var_{a,\l}}_{\infty,\O}\,\leq\,O\Bigl(\displaystyle\frac{1}{\l^{n-4\over 2}}\Bigr)$
and this achieves the proof of the lemma.
\end{pf}
\vfill\eject
\begin{lem}\label{l:ineq}{\rm \cite[lemma 7]{ba2}}
Let $ q > 2$ be given. 
\begin{description}
\item(i) There exists $ \gamma > 1$ such that  for any $ (a_1, \cdots,a_p)\in \bigl(]0,\,+\infty [\bigr)^p$, we have
%\be\label{usefor1}
$$\left( \sum_{i = 1}^p a_i \right)^q \, \geq \, \sum_{i = 1}^p a_i^p \, + \, \frac{\g q}{2} \, \sum_{i \not= j} a_i^{q - 1} a_j
$$
%\ee
\item(ii) There exist $M>0$ such that for any $a$, $b$ in $\R$, we have
$$
%\begin{equation}\label{usefor2}
\bigl|(a\,+\,b)^q\,-\,a^q\,-\,qa^{q-1}b\bigr|\,\leq\,M\bigl(|b|^q\,+\,|a|^{q-2}\min(|a|^2,|b|^2)\bigr).
$$
%\end{equation}
\end{description}
\end{lem}
\begin{lem}\label{usefull}
We have the following formula
\begin{description}
\item(i) $\displaystyle\int_\O\delta_{a,\l}^{2n\over n-4}\,dx\,=\,S^{n-4\over
4}\,+\,O\bigl({1\over\l^n}\bigr)$
\item(ii) $\displaystyle\int_\O\delta_{a,\l}^{n+4\over
n-4}\,\var_{b,\l}\,dx\,=\,\frac{c_2c_n^{n+4\over
n-4}}{\l^{n-4}}H(a,b)\,+\,O\bigl({1\over\l^{n-2}}\bigr)\qquad\mbox{for any
}a,\,b\in\O$ 
 \end{description} 
\end{lem}
\vskip .5truecm
\centerline{\bf Acknowledgements:}
\vskip .1truecm
\noindent The authors would like to thank S.I.S.S.A for its support through postdoctoral fellowships.
 
\end{document}